\documentclass{article}

\usepackage{amsmath,amsfonts,amssymb,amscd,amsthm}
\usepackage{amssymb,amsmath,hyperref}
\usepackage{latexsym}
\usepackage{xypic}
\usepackage{xcolor}
\usepackage{capt-of}
\usepackage{tikz}
\usepackage{pdfpages}
\usetikzlibrary{cd}
\usetikzlibrary{matrix}
\renewcommand{\P}{\ensuremath{\mathbb{P}}}
\newcommand{\Q}{\ensuremath{\mathbb{Q}}}
\newcommand{\p}{\ensuremath{\mathfrak{p}}}
\newcommand{\A}{\ensuremath{{\cal A}}}
\newcommand{\Z}{\ensuremath{\mathbb{Z}}}
\newcommand{\R}{\ensuremath{\mathbb{R}}}
\newcommand{\F}{\ensuremath{\mathbb{F}}}

\renewcommand{\O}{\ensuremath{{\cal O}}}
\renewcommand{\S}{\ensuremath{{\cal S}}}

\DeclareMathOperator{\Aut}{Aut}
\DeclareMathOperator{\Kum}{Kum}
\DeclareMathOperator{\End}{End}
\DeclareMathOperator{\Spec}{Spec}
\DeclareMathOperator{\Stab}{Stab}

\DeclareMathOperator{\Pic}{Pic}

\DeclareMathOperator{\SL}{SL}
\DeclareMathOperator{\PSL}{PSL}

\newtheorem{prop}{Proposition}
\newtheorem{thm}{Theorem}
\newtheorem{lem}{Lemma}

\newtheorem{conj}{Conjecture}
\newtheorem{quest}{Question}

\newtheorem{cor}{Corollary}
\theoremstyle{definition}
\newtheorem{heuristic}{Heuristic}
\newtheorem{expl}{Example}
\newtheorem{defn}{Definition}
\newtheorem*{nota}{Notation}

\newtheorem{rmk}{Remark}
\theoremstyle{remark}
\newtheorem*{ack}{Acknowledgments}
\def\la{\lambda}
\def\eps{\varepsilon}
\def\tK{{\widetilde K}}

\def\Azero{{A_0}} %label for subgroup of Aut coming from char zero
\def\AI{{\overline{A}_0}} %image of Azero in Symmetric group

\newcommand{\xyzw}{[x:y:z:w]}
\newcommand{\sgp}[1]{\ensuremath{\langle#1\rangle}}
\newcommand{\GL}{\operatorname{GL}}
\newcommand{\PGL}{\operatorname{PGL}}
\title{The automorphisms of certain Kummer surfaces act on the points mod $p$
	by the alternating group}
\author{Adam Logan$^1$ $^2$ $^3$ $^4$ 
\and Owen Patashnick$^5$ $^6$ $^7$ $^8$
}
\date{%  
$^1$ School of Mathematics and Statistics, Carleton University, %\\ 4302 Herzberg Laboratories,	1125 Colonel By Drive, Ottawa, ON K1S 5B6, 
Canada\\%
$^2$ ICERM, Brown University, Providence, USA \\% 
$^3$ The Tutte Institute for Mathematics and Computation, %P.O. Box 9703, Terminal, Ottawa, ON K1G 3Z4, 
Canada\\%
$^4$ Email address: \texttt{adam.m.logan@gmail.com}\\ [2ex]%
$^5$ School of Mathematics, University of Bristol, Bristol, %BS8 1TW, 
UK\\%
$^6$ Department of Mathematics, King’s College London, 
UK\\%
$^7$The Heilbronn Institute for Mathematical Research, Bristol, UK\\%
$^8$ Email address: \texttt{owen.patashnick@kcl.ac.uk}\\[2ex]%
\today}
\begin{document}
\maketitle

%\begin{abstract}
{\bf Abstract.}
For the family of Kummer surfaces of the square of an elliptic curve over
the prime field $\F_p$ with $p$ odd, we show that if the automorphism group
acts transitively on the set of points then the action is at least alternating.
This is analogous to a result of Meiri, Puder, and Chen \cite[theorem 1.3]{mpc}
on the points mod $p$ of the Markoff surface.
%\end{abstract}

\section{Introduction}
The classical {\em Markoff problem} is to describe the integer solutions
to the equation 
\begin{equation}\label{eqn:markoff} x^2 + y^2 + z^2 = 3xyz. \end{equation}
%% \footnote{or alternatively to the
%% equation $x^2 + y^2 + z^2 = xyz$ which can easily be reduced to the
%% above equation.}
Given a solution $(a,b,c)$, we obtain a new solution
$(3bc-a,b,c)$ by basic properties of quadratic polynomials,
and similarly for the other two variables.  It is also clear that
negating two of the variables of a solution gives another one, and that
in every solution an even number of the variables are negative.  Accordingly
it is enough to describe the positive solutions.  In some literature the
equation is given as $x^2 + y^2 + z^2 = xyz$; a simple $3$-adic argument
shows that if $(x,y,z)$ is an integral solution to this equation, then
$3|x, y, z$, and so $(x/3,y/3,z/3)$ is an integral solution to 
(\ref{eqn:markoff}).

\begin{defn} A {\em Markoff number} is a positive integer $n$ for which the
(\ref{eqn:markoff}) can be solved with $x = n$.
\end{defn}

We define a graph $G_M$ whose vertices are the Markoff numbers as follows.  If
$(a,b,c)$ and $(3bc-a,b,c)$ are both solutions to the equation, then introduce
an edge from $a$ to $3bc-a$; similarly for the other variables.  
It is known that every positive solution can be reduced to $(1,1,1)$ by 
applying these transformations; in other words, the graph is connected.
Nevertheless, the {\em unicity conjecture} made by Frobenius \cite{frob} 
over a century ago remains intractable, though
there are partial results in \cite[Theorem 0.1]{baragar} and elsewhere:

\begin{conj}\cite{frob} The graph $G_M$ is a tree.
\end{conj}

Recent work, especially \cite{bgs}, \cite{mpc},
has greatly increased our understanding of the
structure of the solutions mod $p$.  Let
$M_p$ be the set of solutions to the Markoff equation mod $p$, excluding
the singular point $(0,0,0)$.  
Let $\Gamma$ be the group generated by the automorphisms
that fix two of the variables as above; let $V$ be the group generated by
the automorphisms that negate two of the variables, and let 
$\Gamma' = \langle \Gamma,V \rangle$
(we have $[\Gamma':\Gamma] = 4$).  We then have the following results:

\begin{thm}\label{thm:big-orbit} \cite[Theorem 1]{bgs}
Fix $\eps > 0$.  Then there exists
$p_\eps$ such that for all primes $p > p_\eps$ there is a $\Gamma$-orbit on
$M_p$ whose complement has order at most $p^\eps$.
\end{thm}

\begin{thm}\label{thm:almost-transitive} \cite[Theorem 2]{bgs}
For all $\eps > 0$, the number of primes less than $n$ for which the action of 
$\Gamma$ on $M_p$ is not transitive is $O(n^\eps)$.
\end{thm}

\begin{thm}\label{thm:alternating} \cite[Theorem 1.3]{mpc} Let 
$p \equiv 1 \bmod 4$ be prime and consider the action of $\Gamma$ on 
$M_p/V$.  Let $G$ be the image of $\Gamma$ in the symmetric
group on $M_p/V$.  If $G$ acts transitively, then it contains
the alternating group.
\end{thm}

The Markoff surface is of degree $3$; accordingly it can be thought of as
a log K3 surface \cite[Chapter 2]{matsuki}.  Integral points on such a surface are
expected to behave analogously to rational points on K3 surfaces.  In addition,
the automorphisms $(a,b,c) \to (bc-a,b,c)$, etc., preserve conic
fibrations and should be thought of as negation automorphisms of
log elliptic fibrations on these surfaces.
It is therefore of interest to investigate whether results such as
Theorems~\ref{thm:big-orbit}, \ref{thm:almost-transitive}, \ref{thm:alternating}
have analogues for certain families of K3 surfaces.

\begin{nota}
Let $X$ be a K3 surface over a number field $L$ and let $\p$ be a prime of
good reduction for $X$ with residue field $\F_\p$.  
Let $X_\p$ be the reduction of $X \bmod \p$
(more precisely, choose a model ${\mathcal X}$ of $X$ over an open subset
of $\Spec \O_L$ containing $\p$ and let $X_\p$ be the special fibre of this
model at $\p$).  There is a map $\Aut X \to \Aut X_\p$ which is usually 
injective, and we may consider the action of $\Aut X$ on the set of
$X(\F_\p)$ of $\F_\p$-rational points $X_\p$.
In this paper we are concerned with the case $L = \Q$ and so when stating results and proofs we will
write $p$ instead of $\p$ and $X(\F_p)$ instead of $X(\F_\p)$.  Note that, although we expect that our results and proofs of same will extend to more general number fields for some choice of K3s, we do not expect them to extend for the K3 surfaces we study in this paper.  
%\owen{OK How's that?} %\owen{Adam, are you ok with the last two sentences?  This is really more of a remark than a notational paragraph I think.  Should we call it a remark instead?}
\end{nota}
%\owen{I left the $X(\F_\p)$'s and $\p$'s that follow as I think you meant them to be there.  Or did you mean for them to all be changed to $X(\F_\p)$s and $p$'s as per the notational remark?}

In analogy
with the results on the Markoff surface above, one can ask about the behaviour
of this action as $\p$ varies, or as $X$ varies within the moduli space of
K3 surfaces with a given lattice polarization.  (A lattice polarization on
a K3 surface $X$ is a fixed embedding of a lattice $L$ into $\Pic X$.
If $L$ is a lattice of rank $1 \le r \le 20$, there is a 
moduli space of K3 surfaces lattice-polarized by $L$ whose dimension
is $20-r$: for the construction and some basic properties see \cite[Section 3]{h-t}.)

In analogy with Theorems \ref{thm:big-orbit}, \ref{thm:almost-transitive},
\ref{thm:alternating}, we might ask the following questions:
\begin{enumerate}
\item Is the action of $\Aut X$ on $X(\F_\p)$ transitive for almost all $\p$?  
\item Is it primitive?
\item If so, can we characterize the permutation group on $X(\F_\p)$ induced by
the action of $X$, and in particular is it the alternating or the symmetric group?
\end{enumerate}

In principle, these questions could be asked about any family of K3 surfaces.
To make them reasonable in practice, the family must have certain properties.  
First, we need a good understanding of $\Aut X$.
This group has been determined by means of a method of Borcherds 
\cite{borcherds} in many cases
(in addition to \cite{kk}, the results of which are fundamental to this paper,
we might mention \cite{kondo} and \cite{ujikawa} among others)
where $\Pic X$ has a complement in the even unimodular lattice $II_{1,25}$
which is a root lattice.  The method is based on the
theory due to Conway and Sloane \cite[Chapter 27]{splag} of the automorphism
group of $II_{1,25}$.
An algorithm of Shimada \cite[Algorithm 6.1]{shimada} can be applied to determine the automorphism
group under much weaker hypotheses, but it is much more difficult to implement.

For certain families of $X$, it is easy to see that some of these questions have
a negative answer.

\begin{expl}\label{ex:solvable-by-finite}
If $\Aut X$ is finite, it certainly cannot act transitively
on arbitrarily large finite sets such as the $X(\F_\p)$.  K3 surfaces with
finite automorphism group have been classified by Nikulin and Vinberg:
see \cite[Theorem 1.10]{acl} for references to complete lists of the possible
Picard lattices.
Even if $\Aut X$ is infinite, these properties may fail for simple reasons.
For example, if $\Aut X$ is virtually solvable, then it cannot
act by the alternating group on $X(\F_\p)$, since the index of the largest
solvable subgroup in $\A_n$ goes to $\infty$ with $n$.  Some families of
elliptic K3 surfaces of Picard rank $3$ with this property
are described by Shimada in \cite[Section 9]{shimada}, \cite{shimada-web}.
\end{expl}

\begin{rmk}\label{rem:tits-alternative}
  Since $\Aut K$ can be embedded into $\GL_n(\Q)$ for suitable $n$,
  %\owen{Is this true for $X$ over $L\ne \Q$?}
   the
  hypothesis that $\Aut(K)$ is virtually solvable is equivalent by a 
  well-known result \cite[Theorem 1]{tits} to the
  statement that $\Aut K$ does not contain a free subgroup of rank $2$.
\end{rmk}

In cases where the automorphism group is large, it may nevertheless fix
a finite set of curves, which prevents any such result from being proved.

\begin{expl}\label{ex:pic-product}
Let $E, E'$ be two elliptic curves that are not isogenous
and do not have complex multiplication,
and let $K$ be the minimal desingularization of
$(E \times E')/\pm 1$.  On $K$ there are eight rational curves
$C_1, \dots, C_4, C_1', \dots, C_4'$
that are fibres of the constant fibrations to $\P^1$ with $j$-invariants
$j(E), j(E')$.  Let $C = \{C_1, \dots, C_4\}$ and $C' = \{C_1', \dots, C_4'\}$.
One checks that every generator $\alpha$ of $\Aut K$, as given in 
\cite[Theorem 5.3]{kk}, preserves the set $C \cup C'$.
More precisely, the action of $\Aut K$ on this set of eight curves is through
the subgroup of $\A_8$ that preserves the partition into two sets of four.
Thus the action of $\Aut K$ on $K(\F_p)$ has an orbit of order at most $8(p+1)$ and
hence is not transitive.
\end{expl}

A well-known result of Babai \cite[Theorem 1.2]{babai} asserts that two random 
elements of $\S_n$ generate at least $\A_n$ with high probability.  Thus one
would expect heuristically that $\Aut X$ should almost always
act on $X(\F_\p)$ at least through the alternating group unless there is a good
reason for it not to.  However,
in order to prove a general statement,
we need some concrete properties of the induced permutations.
For the Markoff surface, this
is done in some cases by means of Jordan's theorem on primitive
permutation groups \cite[Section 8B]{isaacs}.
To apply this result for a general K3 surface, we would need an automorphism
that acts on the $\F_p$-points by a $q$-cycle, where $q$ is prime.  One might
hope to construct such an automorphism by finding an elliptic fibration on the
surface and a prime number $q$ such that exactly one fibre has order $q$ 
(and no other fibre has order a multiple of $q$, though this is automatic from
the Weil bounds if $q$ is not very small).  Then if there is a section passing
through a nonzero point on this fibre, a suitable power of the translation map
will have the desired property.  This is a condition that is easy to check in
practice but not to prove in general.

\begin{expl}\label{ex:tri-inv}
Let $K$ be a {\em tri-involutive K3 surface} over $\F_p$, i.e., a hypersurface of
tridegree $(2,2,2)$ in $(\P^1)^3$.  
Like the Markoff surface, such a surface admits three automorphisms
of order $2$ with (generically) no relations among them.  By Babai's
result, three random involutions in $\S_n$ almost certainly generate at least
$\A_n$; likewise, experiments indicate that this holds
for the action of the group generated by the involutions on $K(\F_p)$ when
$K$ is chosen randomly.
We do not know how to prove any theorem in this direction, and results of
O'Dorney \cite[Theorem 2]{od} show that it fails for special members of the family.
\end{expl}

There are only a few families of K3 surfaces for which enough information
is available to prove a theorem analogous to \ref{thm:alternating}.
We now describe one of these and then introduce some notation.

\begin{defn} For $E$ an elliptic curve over a field, let $\Kum(E \times E)$ be the
minimal desingularization of $(E \times E)/\pm 1$.  This is known as the
{\em Kummer surface} of $E \times E$.  When we discuss properties that depend only
on the Picard lattice and not on the specific $E$, we use the notation
$\tK$.
\end{defn}

\begin{nota}
Let $\Azero$ denote the subgroup of $\Aut K$ that comes from reduction
from $\Kum(E' \times E')$, where $E'$ is a lift of $E$ to $\Q$ without
complex multiplication.  (The automorphism group of $\Kum(E' \times E')$
was determined in \cite[Theorem 5.4]{kk}.
Let $\S_n$ denote the symmetric group on $n$
letters, and let $\S_{X(k)}=\S_{|X(k)|}$.  Similarly $\A_n$ refers to
the alternating group.  An element or a subgroup of $\S_n$ is {\em even} if 
it is contained in $\A_n$.
\end{nota}

Our analogue of Theorem \ref{thm:alternating}, and the main result of
this paper, is the following:
\begin{thm}\label{thm:alternating-kummer} Let $E$ be an elliptic curve
over $\F_p$ with full level-$2$ structure,
where $p>2$.  Let $\AI$ be the image of $\Azero$ in $\S_{K(\F_p)}$.  If $\AI$ is
transitive, then it is $\A_{K(\F_p)}$.
\end{thm}

We have not yet tried to prove results analogous to 
Theorems~\ref{thm:big-orbit}, \ref{thm:almost-transitive} in this setting.
It would be interesting to do so.

The following is an immediate corollary:
%of our theorem (Theorem \ref{thm:alternating}): 
\begin{cor}\label{cor:all-auts}
  Under the assumptions of Theorem~\ref{thm:alternating-kummer}, the image
of $\Aut K$ in $\S_{K(\F_p)}$ contains $\A_{K(\F_p)}$.
\end{cor}

One can of course ask for something stronger:

\begin{quest}
Under the assumptions of Theorem~\ref{thm:alternating-kummer}, 
is the image of $\Aut K$ in $\S_{K(\F_p)}$ equal to $\A_{K(\F_p)}$?
\end{quest}

\begin{rmk}
It seems that the answer to this question is yes for
Kummer surfaces coming from quotients of CM
elliptic curves with small discriminant (e.g. for the cases of CM
discriminant -3 and -4, the automorphism group is generated by maps
that specialize automorphisms of the Kummer of $E \times E$
for generic $E$; cf.~the description of the automorphisms in
\cite[Sections 4.3--4.4]{kk}).  However, we have no reason to
expect any particular answer in general.
\end{rmk}

The following heuristic convinces us that the action of $\Azero$
should always be transitive.
\begin{heuristic}\label{heur:always-transitive}
  First, we will show in Proposition~\ref{prop:trans} that the action
  of $\Azero$ on the $\F_p$-points of a node or line on a suitable
  projective model of $K$ is transitive.
  Given this, let us fix $P \in K(\F_p)$ and an elliptic fibration 
  $\pi: K \to \P^1$ which is the base change
  of a fibration from characteristic $0$ for which the rank of the
  subgroup of the Mordell-Weil group generated by curves in the
  $\Azero$-orbit of a node has rank $\ge 3$ (there are
  fibrations of rank $3$ but none of higher rank).
  With probability $1-O(1/p)$, the fibre containing $P$ will be smooth.
  If the generators can be
  viewed as random points in the fibre, then 
  $P$ is in the subgroup that they generate with probability close to
  $1$, and these events are expected to be independent for different
  choices of $\pi$.  We therefore expect (and experiment agrees) that
  it is easy to find $\pi$ for which $P$ is in this subgroup, and
  therefore that $P$ lies on the reduction of a section in the
  $\Azero$-orbit of a node.
\end{heuristic}
As further support for our belief, we note that we have tested it for
all choices of $E$ for $p < 100$ without finding a counterexample.
In addition, it is supported by the result of Babai cited above.

More generally, we might ask for a description of the circumstances
under which a result such as a weaker version of Theorem~\ref{thm:alternating-kummer} (namely, Theorem~\ref{thm:alternating-kummer-deux}) is
expected to hold.
%\adam{Now that I know about O'Dorney's paper I am less
%confident in this.}

\begin{conj}\label{conj:usually-alternating} The obstructions of Example
  \ref{ex:solvable-by-finite}, Example \ref{ex:pic-product}, and \cite{od}
  are the only meaningful obstructions.
  In other words, if $K$ is a K3 surface over $\F_p$ whose automorphism group 
  is not virtually solvable, has no finite orbits on the set of smooth
  rational curves, and does not fix any nonconstant class in
  $\F_p(K)^*/{(\F_p(K)^*)}^2$, then the image of $\Aut K$ in $\S_{K(\F_p)}$ is
  at least the alternating group for all but finitely many $p$, and
  usually for all $p$.
\end{conj}

Note (cf.~Subsection \ref{subsec:evenness-azero}) that the conditions of
Conjecture \ref{conj:usually-alternating} do not guarantee that the image of $\Aut K$ in
$\S_{K(\F_p)}$ is {\em contained} in the alternating group; that statement
is false even for K3 surfaces which are images of
$\Kum(E \times E)$ by a rational map of degree $2$.

We close the introduction by giving a brief sketch of our proof, comparing
it to that of \cite[Theorem 1.3]{mpc}.  The Kummer surface of $E \times E$
admits elliptic fibrations with constant $j$-invariant, which constitute
a single orbit under the automorphism group.  We view these fibrations as
analogous to the log elliptic fibrations on the Markoff surface above.
The Mordell-Weil rank is $1$ and the torsion is $(\Z/2\Z)^2$.  Let $\tau$
be a generator of the Mordell-Weil group mod torsion, and let $n$ be the LCM
of the orders of the good fibres.  We will show that $\tau^{2n}$ has $16$
cycles of length $p$ in its action on $\Kum(E \times E)(\F_p)$ and fixes all
other points.  This is analogous to the element $\sigma$ constructed in
\cite{mpc} with a single $p$-cycle.

We then have to prove that the automorphism group acts primitively. 
In \cite{mpc} this is straightforward for $p \equiv 1 \bmod 4$ but
requires some estimates of character sums for $p \equiv 3 \bmod 4$.
Our argument uses a nontrivial 
geometric input, namely that there is a pair of disjoint rational curves
$L_1, L_2$ on $\Kum(E \times E)$ such that for all pairs of $\Q$-points 
$(p_1,p_2)$, where $p_i \in L_i(\Q)$, there is a smooth rational curve on
$\Kum(E \times E)$ passing through both.

To finish the proof of the first part of our main result, we need to conclude 
from the primitivity that the action is at least
by the alternating group.  The authors of \cite{mpc} do this for $p \equiv 1 \bmod 4$
by appealing to Jordan's theorem on permutation groups \cite[Section 8B]{isaacs};
for $p \equiv 3 \bmod 4$, they use a deep result of Guralnick and Magaard.
We could perhaps use the same result, but we preferred to use a theorem of
Praeger  \cite[Theorem A]{p} that (unlike the theorem of Guralnick and Magaard) 
does not depend on the classification of finite simple groups.

The second part will be proved by an explicit computation for 
%in which we
%take advantage of the flexibility to
%choose
a set of generators slightly
different from that of \cite{kk}.  The order of
the translation map by an $2$-torsion for an elliptic fibration is
$2$, so we can determine its parity from the number of $\F_p$-points on
$\Kum(E \times E)(\F_p)$ and the number of its $8$ fixed points that
are $\F_p$-rational.  Not all of $\Aut \tK$ is generated by
$2$-torsion translations, however, and so we will have to make
additional arguments on the remaining generators.  To do this, we will
regard them as negations for elliptic fibrations or as automorphisms
that preserve the set of fibres of some elliptic fibration while
acting nontrivially on the base of the fibration.
%without preserving those fibres  individually.

\begin{ack} We thank Michael Fryers for conversations that were essential
in the preparation of this paper and in particular for helping us to
understand \cite{kk}.  We also thank Peter Sarnak, whose inspiring lecture
on the Markoff problem spurred our interest in this subject.
 %\adam{Substantially modified.}
%\both{expand on this and remove the number}
\end{ack}

\section{Geometry and automorphism group of the Kummer surface}
Let $E$ be an elliptic curve over a field $F$ 
of characteristic not equal to $2$ without complex multiplication
(in particular $F$ is infinite).  Our main object of study is the
{\em Kummer surface} of $E \times E$.

\begin{defn} Let $K_0 = K_E$ be the quotient of $E \times E$ by the 
involution $(x,y) \to (-x,-y)$.  Let $\tK$ be its minimal desingularization.
\end{defn}

The surface $K_0$ has ordinary double points ($A_1$ singularities) at the
images of the $16$ fixed points of the involution and is smooth elsewhere.  
Thus $\tK$ is a K3 surface, whose Picard number is $19$.

The Picard lattice of $\tK$ and a set of generators of $\Aut \tK$ are described 
in \cite{kk}.  In particular, the Picard lattice has complement 
$D_4 \oplus A_3$ in $II_{1,25}$; since $II_{1,25} \cong E_8^3 \oplus U$, and
$D_4, A_3$ have complements $D_4, D_5$ respectively in $E_8$, it
follows that $\Pic \tK \cong D_4 \oplus D_5 \oplus E_8 \oplus U$, the generators
being fibral curves and a single section for an elliptic fibration.  This is one
of the most concise descriptions.

\begin{rmk}\label{rem:28-curves}
As in \cite{kk}, there are $28$ distinguished rational curves on the Kummer
surface, which are the images of the curves $x = T, y = T, x-y = T$ on
$E \times E$ for $2T = 0$ and the $16$ exceptional curves.  In
characteristic zero, they arise as
faces of a quasi-fundamental domain for the action of $\O(\Pic \tK)$ on
the positive cone of $\Pic \tK \otimes_\Z \R$.
\end{rmk}

%\adam{This (up to section \ref{subsec:rational-curves}) has been somewhat
%  rewritten.  Please make sure it still makes sense.}\owen{its fine}
The usual quartic model of the Kummer surface \cite[Chapter 3]{cf} degenerates
to a double quadric in the case of a reducible curve, but it can be rescued by
passing to a double cover in $\P(2,1,1,1,1)$, defined by polynomials of
degree $2$ and $4$.  In this model, the $16$ exceptional curves map to
singular points.  For computational purposes, however, it is more convenient
to use a model in $\P^3$ which is analogous to the dual of the Kummer surface
in that the $16$ exceptional curves map to curves of positive degree, while
certain special curves on $E \times E$ are contracted.

We recall some of the classical theory of self-duality
of the Kummer surface \cite[Chapter 4]{cf}.  If $S$ is the Kummer surface of
the Jacobian of a curve of genus $2$, then the linear system of cubics passing
through the $16$ nodes has $4$ sections and defines a map to another surface
with $16$ nodes and $16$ conics, each passing through $6$ of the nodes.
Indeed, if $E_i$ are the classes of the exceptional curves, the linear
system is that associated to the divisor class $3H - \sum_{i=1}^{16} E_i$
and therefore has self-intersection $3^2 \cdot 4 - 16 \cdot 2 = 4$,
intersection $0 - (-2) = 2$ with each of the $E_i$, and
$3 \cdot 4 - 6 \cdot 2 = 0$ with a plane passing through six $E_i$.

For $E \times E$, the usual map to the Kummer degenerates to a double cover
of a quadric in $\P^3$, but the quotient $(E \times E)/\pm 1$ maps naturally
to $\P(2,1,1,1,1)$ with image a surface defined by equations of degree $2, 4$
and having $16$ ordinary double points (the images of the $2$-torsion)
as well as $8$ curves $L_i$ of degree $1$ (the images of $E \times T$ and
$T \times E$, where $T \in E[2]$) and $4$ curves $C_j$
of degree $2$ (the images of the
graphs of translation by $T$ for $T \in E[2]$).
The most obvious analogue of the classical construction, considering forms of
degree $3$ in $\P(2,1,1,1,1)$ vanishing on the $16$ singular points
modulo multiples of the equation of degree $2$, does yield a map to $\P^3$,
but again it has degenerated to a double cover of a quadric.  We now describe
how to construct a 
dual which is a quartic in $\P^3$ with only ordinary double points as
singularities, but the relation between the two models is not symmetrical
as in the classical case.

\begin{prop}\label{prop:dual-product}
Let $D$ be the divisor class
$1/2 (\sum_i L_i + \sum_j C_j) + 1/4 \sum_k E_k$.
The map defined by $D$ is a birational equivalence to a quartic in $\P^3$
taking the $L_i$ and $C_j$ to ordinary double points, which are the
only singularities of the quartic, and the $E_i$ to lines.
\end{prop}

\begin{proof}
The self-intersection number $D^2 = 1/4 \cdot 12 \cdot -2 + 1/16 \cdot 16 \cdot -2 + 1/2 \cdot 12 \cdot 4 \cdot 1 = 4$
is easily calculated.  Since each node is on two
$L_i$ and one $C_j$, we have $(D,E_k) = (2+1) \cdot 1/2 + 1/4 \cdot -2 = 1$;
also each $L$ and $C$ contains four nodes and has self-intersection $-2$ on the
resolution (being a rational curve on a smooth K3 surface), so
$(D,L_i) = (D,C_j) = 0$.  We conclude that $D$ has no base components; since
$D^2 = 4$ and $D$ is effective, this means that $|D|$ has $4$ sections.
It turns out that the
classes of the given curves generate the Picard group of the Kummer surface,
so the sections can be determined and found to be a $4$-dimensional space.
By generic computation, the map defined by $|D|$
is a birational equivalence and the image
has no singularities other than the images of the $L_i$ and $C_j$.
\end{proof}
%\adam{add in sentences to explain relation to duality}
%\adam{This ended up being longer than expected.  Should I cut it down?}\owen{I think it's ok}

Each of the $16$ lines on $K$ passes through $3$ nodes.
Conversely, every set of
$3$ collinear nodes lies on a line contained in $K$ (note that each
such line contains $4$ nodes).  The lines are
the images of the $2$-torsion points of $E \times E$, while the
nodes are the images of the translates of the curves
$E \times 0, 0 \times E$, and the diagonal by $2$-torsion points.
%\end{rmk}

%% Note that the images of the $16$ exceptional curves are of degree $1$,
%% while the other $12$ of the $28$ distinguished curves map to nodes, which
%% are the only singularities.

Perhaps the easiest model to describe is as follows.  More generally, let
$E_i$ for $i \in \{1,2\}$ be elliptic curves over a field $F$
defined by equations of the form
$y_i^2 z_i = f_i(x_i,z_i)$, where the $f_i$ is homogeneous of degree $3$.  Let
$\P(3,1,1,1)$ be weighted projective space with the indicated weights
with the variables labeled $t, u, v, w$.  Then $\Kum(E_1 \times E_2)$ is 
birationally equivalent to the surface in $\P(3,1,1,1)$ defined by
$t^2 = f_1(u,w) f_2(v,w)$, since the latter is the image of
$E_1 \times E_2$ by the map defined by
$(y_1y_2z_1^2z_2^2:x_1z_2:x_2z_1:z_1z_2)$.
In turn this leads to a description of the
elliptic fibrations on $\Kum(E_1 \times E_2)$ with constant $j$-invariant.

\begin{prop}\label{prop:easy-model}
 With notation as in the previous paragraph, the Kummer surface 
$\Kum(E_1 \times E_2)$ admits an elliptic fibration over $\P^1$ whose generic
fibre is isomorphic as an elliptic curve over $F(t)$ to 
$y^2 z = f_1(x,z)$ twisted by $f_2(t,1)$.
\end{prop}

\begin{proof} The desired map is given by the equations $(v:w)$.  The
assertions are easily verified.
\end{proof}

This gives us a convenient way to verify that a surface is birationally
equivalent to $\Kum(E_1 \times E_2)$ or to $\tK = \Kum(E \times E)$.

\begin{prop}\label{prop:quartic-p3} Let $E$ be defined by the Legendre equation
$y^2 = x(x-1)(x-\la)$.  Then the surface $K \subset \P^3$ defined by
\begin{equation}
(\la-1)x^2y^2 - \la x^2z^2 + y^2z^2 + x^2w^2 - \la y^2w^2 + (\la-1)z^2w^2
\end{equation}
has $A_1$ singularities at the coordinate points and the points 
$(\pm 1:\pm 1:\pm 1:\pm 1)$.  It is birationally equivalent to $K_0$
  and $\tK$.
\end{prop}

\begin{proof} The assertion on singularities is easy to check by hand.
For the birational equivalence, it suffices to give an elliptic fibration 
whose generic fibre
is as in Proposition~\ref{prop:easy-model}.  One easily verifies that the
map to $\P^1$ given by $(xy+yz+xw+zw:xy-xz-yw+zw)$ has the desired properties.
\end{proof}

\begin{rmk}\label{rem:right-dual}
  Making a suitable choice of zero section, we find that the rational map
  $\Kum(E \times E) \to K$ of degree $2$ given by the map to
  $t^2 = f_1(u,w) f_2(v,w)$ as above composed with the birational equivalence
  of Proposition~\ref{prop:quartic-p3} contracts the curves
  $E \times T, T \times E, \{(x,x+T)\}$ for $T \in E[2]$ and takes points of
  $(E \times E)[2]$ to lines.  It is thus the same map given in
  Proposition~\ref{prop:dual-product} up to an automorphism of the target
  $\P^3$.
\end{rmk}

For the rest of the paper, when we speak of the Kummer surface, we
will understand either the model $K$ or its minimal
desingularization $\tK$.
Note that it is unnecessary to give an explicit birational map from
$K_0$ to $K$.

\subsection{Rational curves and automorphisms}\label{subsec:rational-curves}
In this section we study the action of $\Aut \tK$ on rational curves and
rational points of $K$.  To simplify the situation we will assume that the
ground field is prime, either $\Q$ or $\F_p$.  As before, we will denote
the ground field by $F$.

In \cite{kk}, Keum and Kondo described (up to a minor error discussed in
Remark \ref{rem:keum-kondo-confused})
 a set of generators for $\Aut \tK$
in the case where $E$ is an elliptic curve over a field of
characteristic $0$ without complex multiplication.  In the case where
$E$ has complex multiplication, and in particular when $E$ is defined
over a finite field, there are additional generators of the
automorphism group, determined by $\End E$ with its Galois action.
These are
determined in \cite{kk} in the case of the endomorphism ring having
discriminant $-3$ or $-4$, but it seems very difficult to give any
kind of general description.

\begin{defn}\label{defn:AutKChar0}
  Let $\Azero$ denote the image in $\Aut \tK$ of the group of
  automorphisms described in \cite[Section 4.2]{kk} (Case II in the
  notation of their paper).
  \end{defn}

The rational curves in the $\Azero$-orbit of a node are of special
importance for this paper.

\begin{defn}\label{def:basic-curves}
  Let $C$ be a smooth rational curve on $\tK$.  If $C$ is in the
  $\Azero$-orbit of a curve lying above the line on $K$ which
  is the image of the origin of $E \times E$,
  we say that $C$ is {\em basic}. %\adam{I changed this name from ``knotty'',which I didn't feel like explaining when I gave the talks.  We can discuss.}
\end{defn}

\begin{expl} All $16$ lines on $\tK$ are basic.
  Indeed, translation by any $2$-torsion point of $E \times E$
  descends to an automorphism of $\tK$, and the group of these
  automorphisms acts simply transitively on the set of these $16$
  curves.  Later we will see that the $12$ exceptional curves above
  the nodes of $K$ are basic as well.  In characteristic $0$, in fact,
all smooth rational curves on $\tK$ are basic, but we will not need this.
\end{expl}

\subsubsection{A concise set of generators of
  $\Aut \tK$}\label{alternategens}

The set of generators for $\Aut \tK$ given in \cite{kk} is closely
related to the embedding of the Picard lattice into $II_{1,25}$ and is
therefore useful for theoretical purposes.  However, the number
of bounding hyperplanes in this case is $119$ \cite[Lemma 3.6]{kk},
and if we took a generator for each of these in addition to generators
for the finite group $G_0$ acting trivially on the roots
\cite[Lemma 3.5]{kk}, our generating set would be inconveniently large.
Here we discuss a set of generators for $\Aut \tK$ that is smaller
and quite convenient for computation.  Essentially our generators consist
of generators for their finite group $G_0$ and representatives of every
other type of generator given in \cite{kk} up to conjugacy by $G_0$.
We found it necessary to modify one of Keum and Kondo's generators for reasons
discussed in Remark~\ref{rem:keum-kondo-confused}.
We start by recalling the notation of \cite{kk} for the nodes
and lines of $K$.

\begin{defn}\label{def:node-line-names}  Let $O = T_0$ be the origin
  of $E$ and let $T_1, T_2, T_3$ be the points of order $2$ (recall
  that we are assuming that these are rational).  Specifically, we
  take $T_1, T_2, T_3$ to be the points $(0:0:1), (1:0:1), (s:0:1)$
  respectively. Define $G_{ij}$ for
  $0 \le i,j \le 3$ to be the line on $K$ which is the image of the
  point $(T_i, T_j) \in E \times E$.  Let $E_i, F_j$ be the nodes
  which are the images of $E \times T_i$ and $T_j \times E$
  respectively, and let $D_i$ be the translate of the diagonal of
  $E \times E$ by $T_i$.
\end{defn}

\begin{defn}\label{def:constant-fibs} The {\em horizontal}, {\em
    vertical}, and {\em diagonal fibrations} are the elliptic
  fibrations induced on $\tK$ by the maps $E \times E \to E$ taking
  $(P,Q)$ to $Q, P, P-Q$ respectively.  These are referred to
  collectively as the {\em constant fibrations}, because the fibres have
  constant $j$-invariant, although there are other fibrations with the
  same property.  Each of them has four fibres of type $I_0^*$.
\end{defn}

\begin{rmk}\label{rem:pic-tk}
  If $E$ has no complex multiplication, then the Picard group of $\tK$ is generated by the
  classes of lines and nodes.  The relations are those induced by the
  equalities of the classes of reducible fibres of the constant
  fibrations.  For example, two of the $I_0^*$ fibres of the diagonal
  fibration are supported on $\{D_0, G_{00},G_{11},G_{22},G_{33}\}$
  and $\{D_1, G_{01},G_{10},G_{23},G_{32}\}$.  We thus obtain the
  relation
  $2[D_0]+[G_{00}]+[G_{11}]+[G_{22}]+[G_{33}] =
  2[D_1]+[G_{01}]+[G_{10}]+[G_{23}]+[G_{32}]$,
  where $[C]$ denotes the Picard class of the curve $C$.  In this way
  we obtain $9$ relations; these are independent, so they cut down the
  N\'eron-Severi rank of $\tK$ to~$19$.
\end{rmk}

We start by describing the simpler generators on our list.  We relate each of
these to an automorphism of $E \times E$, give polynomials defining
them, and describe them in terms of elliptic fibrations on the surface.

\begin{defn}\label{def:a-to-e}
  The automorphisms $a, b, c, d, e$ are induced by automorphisms of
  $E \times E$ as follows:
\begin{itemize}
    \item $a$ is induced by the map $(P,Q)\mapsto(Q,P)$, 
    \item $b$ by $(P,Q)\longmapsto(P,-P-Q)$, 
    \item $c$ by $(P,Q)\longmapsto(P,Q+(0,0))$,
    \item $d$ by $(P,Q)\longmapsto(P,Q+(1,0))$, and    
    \item $e$ by $(P,Q)\longmapsto(-P,Q)$.
\end{itemize}
%Recall that $(0,0)$, $(1,0)$ and~$(\lambda,0)$ are the points of 
%order~$2$ on~\El.  
\end{defn}

From this it is possible to see that the subgroup $\sgp{a,b,c,d,e}$ is
isomorphic to $\PGL_2(\Z)\ltimes C_2^4$, where, on $E\times E$, each generator of $C_2^4$ acts by translation by a two-torsion point, $\PGL_2(\Z)$ acts in the obvious way, and $\PGL_2(\Z)$ acts on $C_2^4$ by its action on the two-torsion points.
%\owen{I didn't put in the comment about $\PGL_2(\Z)\mapsto \PGL_2(\F_2)$ as I think it was already clear; see what you think}
%\adam{What I was objecting to was that $\S_3$ doesn't act on $C_2^4$ in an
%  obvious way.  I rewrote the next paragraph a bit; does it look reasonable
%  to you?}\owen{yeah that looks good}

\begin{defn}\label{def:h} Let $H$ be the subgroup $\sgp{a,b,c,d}$.
\end{defn}

%%% b is what we usually call ab
Note that $a, b$ generate a subgroup of $\PGL_2(\Z)$ isomorphic to $\S_3$.
It is immediate to verify that $H$
has order~$96$ and is isomorphic to $\S_3\ltimes C_2^4$,
where $\S_3$ acts on $C_2^4$ by restricting the action on $\PGL_2(\Z)$; 
it consists of all those automorphisms which are restrictions
of linear maps $\P^3\to\P^3$.  Thus as maps $K \to K$:
\begin{align*}
    a:\xyzw&
    \longmapsto
    [x+y+z+w:x+y-z-w:x-y+z-w:x-y-z+w],\\
    b:\xyzw&
    \longmapsto
    [-x:y:z:w],\\
    c:\xyzw&
    \longmapsto
    [-x:-y:z:w],\\
    d:\xyzw&
    \longmapsto
    [-x:y:-z:w].
\end{align*}

The formula for $e$ is
\begin{align*}
    e:\xyzw\longmapsto
    &[x(-x^2+y^2+z^2+w^2)+2yzw:\\
    &\;\phantom{{}-{}}y(x^2-y^2+z^2+w^2)+2xzw:\\
    &\;\;\phantom{{}-{}}z(x^2+y^2-z^2+w^2)+2xyw:\\
    &\;\;\;\phantom{{}-{}}w(x^2+y^2+z^2-w^2)+2xyz].
\end{align*}

In terms of \cite[Lemma 3.5, Case II]{kk}, our $a, b, c, d$ generate a
subgroup of $G$ that contains $G_0$ with index $2$.  On the other
hand, $e$ is not defined by linear equations and does not belong to
$G$, because it takes the image of the diagonal on $E \times E$ to the
image of the antidiagonal.

However, we can describe $e$ in terms of the vertical fibration, which
it obviously preserves.  Since $e$ preserves the curves $E_i$, which
are sections of the vertical fibration, and it is not the identity, it
must be the negation with respect to one of the sections $E_i$ as zero
section (since the difference of two of these is $2$-torsion the
negation maps for all of them are equal).  The map $e$ is
the $f_{r'}$ constructed in \cite[Case II-(1)]{kk}.  Since the action of
$H$ on Leech roots (equivalently, bounding hyperplanes)
defining a $D_4 \oplus D_4$ extension is transitive,
only one such generator is needed.

%% Moving on to $f$, we consider the elliptic fibration $\rho$ defined by
%% $(x+z:y+z)$.  Let $S_0, S_1$ be the (rational curve above the) singular point 
%% $(0:0:0:1)$ and the line $y + w = z + w = 0$; both are sections
%% of this fibration.  Let $f$ be the translation automorphism for $\rho$ that 
%% takes $S_0$ to $S_1$.  It is defined by the equations
%% $$\begin{aligned}
%%     &x^2y - xy^2 - x^2z - xz^2 - y^2w + 
%%         z^2w + 2xw^2 - yw^2 + zw^2,\\
%%     &-x^2y + xy^2 + x^2z - 2y^2z - xz^2 - 
%%         y^2w + z^2w + yw^2 + zw^2,\\
%%     &x^2y + xy^2 - x^2z - xz^2 - 2yz^2 + 
%%         y^2w - z^2w + yw^2 + zw^2,\\
%%     &x^2y + xy^2 + x^2z - xz^2 - 2x^2w + 
%%         y^2w + z^2w - yw^2 - zw^2.\\
%%   \end{aligned}$$

%% \begin{rmk}\label{rem:same-or-opposite-parity}
%%   The element $f$ is closely related to the automorphism $f_{r'}$
%%   constructed in \cite[4.2, Case II-(2)]{kk}.  A minor change of coordinates
%%   matches $\pi$ with $\rho$ and identifies $E_1, G_{43}$ as defined in 
%%   \cite[Sections 2]{kk} with $S_0, S_1$ as introduced here.  However,
%%   the $f_{r'}$ of \cite{kk} is a translation of infinite order, whereas
%%   our $f$ is a negation.

%%   The importance of $f$ arises in part from the fact that it changes the
%%   parity of the degree of all smooth rational curves, unlike our
%%   other generators, which preserve it.  (This implies that the action of
%%   $\Aut \tK$ on the set of smooth rational curves is not primitive.)
%% \end{rmk}

Next we consider the elliptic fibration $\rho$ defined by
$(y-w:z-w)$.  For a suitable map from $E \times E$ to $\tK$,
the coordinates of $E_1$ are $(1:0:0:0)$ and $G_{43}$ is defined by
$x-z = y+z = 0$.  The map called $f_{r'}$ in \cite[Section 4.2, Case II-(2)]{kk},
and which we denote by $f$,
is the translation for this fibration that takes $E_1$ to $G_{43}$.
In terms of coordinates it is defined by
$$\begin{aligned}
x^2y + xy^2 - x^2z + xz^2 - y^2w + z^2w - 
        2xw^2 - yw^2 + zw^2,\\
    -x^2y - xy^2 + x^2z - 2y^2z + xz^2 - y^2w + 
        z^2w + yw^2 + zw^2,\\
    x^2y - xy^2 - x^2z + xz^2 - 2yz^2 + y^2w - 
        z^2w + yw^2 + zw^2,\\
    x^2y - xy^2 + x^2z + xz^2 - 2x^2w + y^2w + 
        z^2w - yw^2 - zw^2\\
\end{aligned}$$
Again the action on Leech roots defining a $D_4 \oplus A_4$ extension is transitive
and we do not need any more generators of this type.
%% In terms of our usual generators this would be called afab^2.
%% So we have to change f to afba everywhere it appears later.
%% Then we have to change all the b's to ab's, so in total f becomes afaba.

We continue with descriptions of $g, h, i$.  It no longer suffices to
choose a single representative map corresponding to $f_r'$ as defined
in \cite[Section 4.2, case II-(3)]{kk}.  Rather, we need three,
one for each of the orbits of $H$ 
on the set of Leech roots corresponding
to root systems of type $D_5 \oplus A_3$.  We have made a somewhat arbitrary
choice, seeing no particular reason to prefer one set of orbit representatives
to another.

Like the $f_{r'}$ of \cite[Section 4.2, case II-(3)]{kk}, 
our maps are related to fibrations with
reducible fibre types $I_2^*,I_2^*,I_4,I_2,I_2$ and full $2$-torsion.
Instead of a translation map, however, we will use a negation.
In particular, consider the $\tilde D_6$-configuration supported on
the resolutions of $(-1:-1:-1:1)$ and $(1:-1:-1:1)$ and on the lines
$$\begin{aligned}
x-w=z+w=0,\quad &x+w=y+w=0,\quad y+w=z+w=0,\\
&\quad x+w=z+w=0,\quad x-w=y+w=0.\\
\end{aligned}
$$
Two of the sections are the nodes $N = (-1:-1:1:1)$ and $N' =(-1:1:-1:1)$, 
and we let $g$ be the negation that interchanges these two: in other words,
the automorphism that on every nonsingular fibre takes a point $P$ to
$N+N'-P$.  The fibration may be
defined by $(yz+w^2:(y+z)w)$.

\begin{rmk}\label{rem:keum-kondo-confused}
  We modify the set of generators given in \cite{kk} because we do not 
understand the proof given by Keum and Kondo that their generators 
of $\Aut \tK$ actually
generate this group.   The element $f_{r'}$ that they construct in Case II-(3)
does not map the quasi-fundamental
domain to an adjacent isomorphic copy, but rather to a copy at distance $2$.
On the other hand, the automorphism $g$ takes the quasi-fundamental
domain to an isomorphic copy that shares a face of codimension $1$, and
the proof of \cite[Theorem 5.4]{kk} is valid if $f_{r'}$ is replaced by $g$.

% note again that their b is our usual ab
The set of generators given in \cite[Case II]{kk} do in fact generate 
$\Aut \tK$; indeed, their $f_{r'}$ in case II-(3) is our $baeabg$, so it
may be used as a generator in place of $g$.
We thank Michael Fryers for very illuminating discussions of these points.
\end{rmk}

With this done, the definitions of $h, i$ are similar; we just need to
permute $y, z, w$, which corresponds to relabeling the torsion points
of $E$.  That is, for $h$, we consider
the fibration defined by $(zw+y^2:(z+w)y)$ and the negation that exchanges
$(-1:1:1:-1)$ with $(-1:1:-1:1)$.  For $i$, we use the fibration
$(yw+z^2:(y+w)z)$ and the nodes $(-1:1:1:-1),(-1:-1:1:1)$.

Our $j, k, \ell$ are automorphisms that preserve fibrations with an
$I_{16}$ fibre as in \cite[Section 4.2, case II-(4)]{kk}, and like them
they are neither negations nor translations, but rather act nontrivially
on the base of the fibration.
 
Again there
are $3$ orbits of bounding hyperplanes corresponding to root systems
of type $D_8$ and so we need three generators for this case.
We add to the description there by giving the reasonably elegant equations
\begin{align}
  \nonumber &(\la x y^3 + \la y^3 z - \la x y z^2 - \la y z^3 + \la x y^2 w + \la y^2 z w - \la x z^2 w - \la z^3 w:\\
  \label{eqn:j} &\quad x y^2 z + (\la+2)/2 y^2 z^2 - (\la/2) z^4 + \la y^3 w + 2 x y z w - (\la-2) y z^2 w + \\
  \nonumber &\qquad(\la/2) y^2 w^2 + x z w^2 - (\la-2)/2 z^2 w^2)\end{align}
for the fibration and identifying some small sections as the lines given
by $$y+w=z-w=0,\quad y+w=z+w=0,\quad x+w=z+w=0,\quad x-w=z-w=0.$$
Then $j$ acts on the base as $(s:t) \to (-s:t)$
(in particular preserving the $I_{16}$ fibre at $(0:1)$)
while also preserving these four sections and fixing every component of the
reducible fibre as a set.

For $k,\ell$ the equations of the fibration are not quite as nice, but
it remains true that the automorphism acts nontrivially on the base,
fixes four lines that are sections and the $16$ vertical curves, and is
determined uniquely by this information.  We relegate the equations to
Appendix \ref{app:eqns}.

\begin{rmk}\label{rem:jkl-big-degree}
In our model the maps $j, k, \ell$ can only be defined by polynomials of
degree at least $28$.  Fortunately, such polynomials are not needed in this paper.
These generators do not appear in the automorphisms that
we will use to study the pointwise stabilizer of a curve in
Proposition \ref{prop:trans-2}, and we can determine the fixed points
of $dj, cdk, c\ell$ without explicit equations, which is enough to prove that
they induce even permutations (Theorem \ref{thm:alternating-kummer-deux}).
\end{rmk}

\subsection{Actions on curves in $\tK$}
\begin{lem}\label{lem:gen-aut}
  Let $C$ be a basic curve.  Then there is a basic curve $C'$ with
  $(C,C') = 1$ such that $\Azero$ is generated by the stabilizer of
  $C$ together with an element $\alpha$ such that $\alpha(C) = C'$.
\end{lem}

%% changed dfd to dafabad, see above
%% again, ab is what we usually call b (and aba here is right)
\begin{proof} We do this computationally.  Let $C = F_0$, let
  $C' = G_{00}$, and take $\alpha = dafabad$.  We find that
  $\{aba, aca, ada,e,aga,h,i,j,k,\ell\} \subset \Stab_{\Azero}(C)$, while
  $a \in \Stab_{\Azero}(C')$.  It follows that all of
  $a, b, \dots, \ell$ belong to $\sgp{\Stab_{\Azero}(C),dafabad}$ except
  possibly for $f$.  But a subgroup containing not only $dafabad$ but also
  $a, b, d$ contains $f$ as well.
%Basically, we are showing that 
%$\Azero$ is generated by the stabilizer of a curve together with
%something that takes the curve to a curve that meets it in one point.  
\end{proof}

\begin{lem}\label{lem:chainbasic}
  Let $C, C'$ be basic curves.  Then there are basic curves $C = C_0,
  C_1$, $\dots$, $C_n = C'$ such that $(C_{i-1},C_i) = 1$
  for $0 < i \le n$.
\end{lem}

\begin{proof}
  It suffices to prove this for $C = F_0$.
  Let $\alpha \in \Azero$ be such that $\alpha(C) = C'$ and let
  $\delta = dafabad$, so $(\delta C,C) = 1$.  Let $D = \delta \Stab_{\Azero}(C)$.
  Using the lemma, we write $\alpha$ as a product $d_n d_{n-1} \dots d_1$
  of elements of $D$.  Now define $C_i = d_n d_{n-1} \dots d_{n-i} C$ for
  $0 \le i \le n$.

  By definition the $C_i$ are basic curves, so it suffices to show
  that $(C_{i-1},C_i) = 1$.  This holds, because $C_{i-1}$ and $C_i$
  are the images of $C$ and $d_{n-i} C$ by the automorphism
  $d_n d_{n-1} \dots d_{n-i+1}$, and $d_{n-i} C$ = $\delta C$ for all $i$.
\end{proof}

%% \begin{lem}\label{lem:chainbasic}
%% Any two basic curves are connected by a sequence of basic curves
%% each of which intersects its neighbor in a rational point.
%% \end{lem}

Perhaps the most obvious automorphisms of $\tK$ are those induced by the action
of $\SL_2(\Z)$ on $E \times E$.  These are already sufficient to prove that
the action of the automorphism group on the $\Q$-rational points of basic
rational curves is transitive.

\begin{prop}\label{prop:trans} Suppose that $E$ is defined over $F = \Q$,
respectively 
$F = \F_p$ with $p$ prime, and let $C$ be the rational curve given by resolving
the singularity of $(E \times E)/\pm 1$ at the origin.  Then the action of
the stabilizer of $C$ in $\Azero$ on the $F$-points of $C$ is transitive, respectively 
$3$-transitive.
\end{prop}

\begin{proof} Indeed, the subgroup $\PSL_2(\Z)$ of the automorphism group
acts by the usual action of $\PSL_2(\Z)$ on $\P^1(\Q)$ or $\P^1(\F_p)$, which
is well-known to be transitive, respectively $3$-transitive.
\end{proof}

\begin{cor}\label{prop:trans-basic} The same statement holds when $C$
  is replaced by any basic curve.
\end{cor}

\begin{proof} Let $C$ be the curve from Proposition \ref{prop:trans}.  Then
  the basic curve in question is $\alpha(C)$ for some $\alpha \in  \Aut \tK$.
  If $\beta$ fixes $C$ and has a given
  action on it, then the action of $\beta^\alpha$ on $\alpha(C)$ is
  the same.
\end{proof}

\begin{rmk}\label{rem:finiteorbitssrc}
By a general result of Sterk \cite[Lemma 2.4 and previous discussion]{sterk}
 there are only finitely many
orbits of smooth rational curves on $\tK$.
\end{rmk}

\begin{rmk}\label{rem:nodesorbit}
  The $12$ nodes are basic.  Indeed, in Section~\ref{alternategens} we
  described an automorphism $f$ and asserted that it changes the parity
  of every smooth rational curve.  More precisely, it takes every node to
  a line; this can be worked out from explicit equations (it is easier to
  show that $12$ lines go to nodes) or by analyzing the fibration.
  Since all lines are basic, the same follows for nodes.
\end{rmk}

\begin{cor}\label{cor:trans-all}
Let $C$ be any of the $28$ rational curves of Remark~\ref{rem:28-curves}.
Then the action of $\Aut \tK$ on $C(F)$ is transitive for $F = \Q$ and
$3$-transitive for $F = \F_p$.
\end{cor}

%% \begin{defn}\label{defn:f} Let $f$ be as defined in the proof of Lemma
%%   \ref{lem:nodesorbit}.
%% \end{defn}

%% \owen{this discussion of $f$ should be moved back to the preceeding
%%   subsubsection \ref{alternategens}}

We are now ready to begin our discussion of the permutation
representation of $\Aut \tK$ on $\tK(\F_p)$.

\begin{lem}\label{lem:Block2points}
Let $C$ be a basic curve, and let $B$ be a block for the action of
$\Aut \tK$ on $\tK(\F_p)$ that contains two points of $C(\F_p)$.
Then $B$ contains all points of $C(\F_p)$.
\end{lem}

\begin{proof}
  Let $P, Q \in C(\F_p)$ belong to $B$, and let $R$ be any other point
  of $C(\F_p)$.  By Corollary~\ref{cor:trans-all} we can find an
  element of $\Aut \tK$ fixing $P$ and taking $Q$ to $R$.  Thus,
  if the action preserves a partition for which $P$ and $Q$ are in the
  same part, then $R$ must be in that part as well.
\end{proof}

Perhaps more surprisingly, we have the following result.
\begin{prop}\label{prop:trans-2}
Let $C_1, C_2$ be rational curves on $\tK$ which are both nodes or both 
lines in our model, and suppose that $F \ne \F_3$.
Then the action of $\Aut \tK$ on $C_1(F) \times C_2(F)$ is
transitive.  In particular, for every pair of rational points $(P_1, P_2)$
in $C_1(F) \times C_2(F)$ there is a basic curve passing through
$P_1, P_2$.
\end{prop}

\begin{rmk} Obviously if $C_1, C_2$ are both nodes they are disjoint.  
In fact the same holds if $C_1, C_2$ are both lines; it can be checked that
the intersection point of any two lines in our model that meet in $\P^3$ is
a node of the surface, and $C_1, C_2$ are disjoint on $\tK$.
\end{rmk}

\begin{proof} 
First, there is a unique $\Aut \tK$-orbit of pairs $\{C_1, C_2\}$ as described.
This can be checked by writing down the action of the generators
% no need to change this
$a, b, \dots, f$ on the curves of degree $0, 1$ and verifying that the graph
whose vertices are the pairs and whose edges join pairs related by a single
generator is connected.   (Note, however, that there are multiple orbits
of pairs of disjoint rational curves without any restriction on the degree.)

Thus it suffices to consider the case where $C_1, C_2$ are the curves in
the resolutions of the nodes at
$(1:0:0:0)$ and $(1:1:1:1)$ respectively.  There exists a pair of
points $(P_1,P_2)$ with $P_i \in C_i$ for $i = 1, 2$ that lie on the
same rational curve: for example, take the $P_i$ to be the intersection
of the exceptional curves with the line $y = z = w$.  We will show that
the pointwise stabilizer of $C_1$ acts transitively on $C_2(F)$.  Since
$aC_1 = C_2$ and $aC_2 = C_1$, it follows by conjugating by $a$ that the
action of $\Stab_{C_1} \cap \Stab_{C_2}$ on pairs $(P_1,P_2)$ is transitive.

In particular, the words 
\begin{align}
\label{eqn:longword}
%&eacdfacdafabdaicbbacgbda, bcgbcdbbdaicbafdbgadfacbacdabacgbcdb,\\
%\nonumber &\quad bbdaicafacdafafdbgadfacicbbafdbgadfd
% changed f to afba then b to ab in these strings, then cancelled things
% (and checked that if we change f to afabb and then b to ab in the
% generators, these do the right thing)
\nonumber &eacdafabcdfbadaicbcgabd, abcgdacdicabfabadabgadafcdagabcdab, \\
&\quad badaicfabcdfabafabadabgadafabcicbafabadabgadafabad
\end{align}
define automorphisms (using Magma's convention of composition from left to
right so that the first of these means that we apply $e$, then $a$, etc.)
that fix $C_1$ pointwise.  On the other hand, let us parametrize $C_2$
so that the intersections with the lines $x=y=z, x=z=w, y=z=w, x=y=w$
correspond to the points $(1-s:1), (1:1), (1:0), (0:1)$ respectively.
Then the action of these automorphisms on $C_2$ is given by the matrices
\begin{equation}
\begin{pmatrix}1&0\\-1&1 \end{pmatrix},\quad
\begin{pmatrix}2&3\\-1&-1 \end{pmatrix},\quad
\begin{pmatrix}0&1\\-s(s-1)^3/(s^2-s+1)&s-1 \end{pmatrix}.
\end{equation}
The first two of these generate the group $\Gamma^0(3)$ of elements
of $\SL_2(\Z)$ whose upper right entry is a multiple of $3$.  This has
two orbits on $\P^1(\Q)$, namely the points that do or do not reduce to
$\infty \bmod 3$ (according to our convention this group acts on the right).
These orbits are exchanged by the third generator, which takes $(1:0)$ to 
$(0:1)$.  It follows that the action of the group is transitive in the case
$F = \Q$.  If $F = \F_p$ and $p \ne 3$, 
the argument is even simpler, since the first 
generator already takes $(a:1)$ to $(a-1:1)$ and the second takes
$(1/3:1)$ to $(1:0)$.
\end{proof}

\begin{rmk}\label{rem:how-to-find}
We pause to comment on how we found and verified the words used in the proof.
First, given an automorphism that fixes a curve (as a composition
of generators given by equations), we would like
to describe its action on that curve.  In favourable situations, where the
curve has dimension $1$ in our model rather than being an exceptional curve
and the automorphisms do not send it through a singularity, we can calculate
the images of points directly.  In order to find the image of a point lying
above a singularity without explicit desingularization, one approach is to
find an elliptic fibration for which the curve in the resolution passing through
the desired point is a section, while all other curves in the same resolution
are fibral (the second condition is vacuous when the singularity is an ordinary
double point, as in this paper).  Automorphisms can then be applied to this 
elliptic curve, and the image will still be of genus $1$ at the end.  One
can keep track of the image of the desired point by mapping an infinitely
near point of sufficiently high order.  There is no difficulty in doing this
because the point near to $n$th order to a smooth point on a curve is unique.
By doing this calculation for $3$ points on the curve we find the action,
and we can verify it by checking additional points.

In order to find the words, we start by the standard computation to write down
the matrix representing the
action of an automorphism given by translation or negation on an elliptic
fibration on the Picard lattice.  Having done that, a simple search finds
a large number of short words that fix the curve as a set.  Computing their
actions on points as above, we find some short products of these words that
act as the identity on $C_1$ and that have simple matrices for $C_2$.
Examining these matrices, we found some that could be used to prove the
desired transitivity.
\end{rmk}

Fix a constant fibration $\pi: \tK \to \P^1$: for concreteness, we choose to
define it by $((x-z)(y-w):(x-y)(z-w))$.  Let $S$
be the exceptional curve above $(1:0:0:0)$.  Let $\tau$ be the
translation automorphism taking $S$ to the exceptional curve above
$(-1:1:1:1)$.  Then it is easily checked that $\tau$ is a generator of
the Mordell-Weil group of $\pi$ modulo torsion, that $\tau^2$ fixes
every component of the bad fibres of $\pi$, and that $(S,\tau^2 S) =
0$ so that $\tau^2$ has no fixed points on good fibres of $\pi$.  The
action on the central component of the reducible fibre fixes the
points of intersection with all $4$ of the other components, so the
central component of the reducible fibre is fixed pointwise.  The non-central
components are
the $16$ lines of $K$.

\begin{lem}\label{lem:tau2action}
The permutation given by $\tau^2$ on the $\F_p$-points of a non-central component of a bad fibre has $1$ cycle of length $p$ and $1$ fixed point.
\end{lem}

\begin{proof}
With the choices made above, we have
$\tau = b^2abe$.  One verifies, for example, that $\tau^2$ restricts
to the automorphism of the line $(-x:x:x:y)$ given by $(-x:x:x:y) \to
(y:-y:-y:x+2y)$, which is equivalent by an easy change of coordinates
to the automorphism of $\P^1$ given by $(x:y) \to (x+y:y)$.  This 
automorphism exhibits the behaviour claimed for it.  The other
$15$ lines reduce to this one because $\tau^2$ commutes with 
$c, d, aca, ada$, and these generate a group that acts transitively on the
lines.
\end{proof}

\begin{lem}\label{lem:good-power-of-tau}
  There is a positive integer $n$ such that $\tau^{2n}$
  has $16$ cycles of length $p$ and fixes all the other points
  in its action on $\tK(\F_p)$.
\end{lem}

\begin{proof}
  The smooth fibres of $\pi$ are elliptic curves over $\F_p$ with full
level-$2$ structure, so the number of points on them is not a multiple
of $p$.  Let $n$ be the LCM of the orders of the smooth fibres.
Since $\tau^2$ acts by translation on a smooth fibre, the orbits are all
of the same size, which divides $n$, so $\tau^{2n}$ fixes every point
on a smooth fibre.  On the singular fibres, we have already seen that
$\tau^2$ fixes the central components and acts on the outer components
with a $p$-cycle and a single fixed point.  The same holds for
$\tau^{2n}$, since $n \nmid p$. 
\end{proof}

\begin{rmk}\label{rem:like-mpc}
  The element $\tau^{2n} \in \Azero$ is 
analogous to the element $\sigma$ constructed by Meiri and Puder in
the proof of \cite[Theorem 3.1]{mpc}.%Meiri and Puder.
\end{rmk}

\begin{lem}\label{lem:two-makes-everything}
  Assume that $\Azero$ acts transitively but not primitively
  on $\tK(\F_p)$; let ${\mathcal B}$ be a partition that exhibits
  the imprimitivity.  Then there is no block of ${\mathcal B}$ that
  contains two points on the same basic curve.
\end{lem}

\begin{proof}
  Let $C_0$ be the basic curve two of whose points are in the same
  block $B_P$, and let $C'$ be any basic curve.  By Lemma
  \ref{lem:chainbasic} there is a chain $C = C_0, C_1, \dots, C_n =
  C'$ of basic curves with $(C_{i-1},C_i) = 1$ for all $i$.  Every
  $\F_p$-point of $C$ is contained in $B_P$, by Lemma~\ref{lem:Block2points},
  so for all $i$, all
  points of $C_i$ are in the same block.  But this block does not
  depend on $i$ because $C_{i-1}(\F_p) \cap C_i(\F_p) \ne \emptyset$.
  Therefore every $\F_p$-point on a basic curve is in $B_P$.  But we
  are assuming that $\Azero$ acts transitively on $\tK(\F_p)$, so this
  means that $B_P = \tK(\F_p)$.
\end{proof}

\begin{lem}\label{lem:primitive}
Assume that $\Azero$ acts transitively on $\tK(\F_p)$.  Then the
permutation action of $\Azero$ on $\tK(\F_p)$ is primitive.
\end{lem}

\begin{proof}
  Suppose that ${\mathcal B}$ is a non-discrete
  partition preserved by the action.
  We must show that ${\mathcal B}$ is the trivial partition.
  Let $P \in \tK(\F_p)$ be a point not fixed by $\tau^{2n}$, let $C$ be the
  line containing $P$, and let $B_P$ be the block containing $P$.
  We may assume that $B_P \ne \{P\}$, so let $Q \ne P \in B_P$.
  If $Q \in C$, then ${\mathcal B}$ is trivial by Lemma
  \ref{lem:two-makes-everything}.  If $Q$ belongs to another line,
  then there is a basic curve joining $P$ and $Q$ by Proposition
  \ref{prop:trans-2} and again the triviality follows.  Otherwise
  $Q$ is fixed by $\tau^{2n}$.  Therefore $\tau^{2n}(B_P)$ is a block
  whose intersection with $B_P$ is nonempty (it contains $Q$),
  so $\tau^{2n}(P) \in B_P$.  But
  $\tau^{2n}(P)\in C$ and $P \in C$, so Lemma \ref{lem:two-makes-everything}
  again applies.
\end{proof}
  
%%   We
%%   consider two cases: either $\tau^2(B_P) = B_P$ or not.

%%   If $\tau^{2n}(B_P) = B_P$, then $B_P$ contains two points of
%%   $C$, so by lemma \ref{lem:tau2action} every
%%   point of $C$ except possibly the point of intersection with the
%%   central component of the fibre containing $C$ is in $B_P$,
%%   and, by lemma \ref{lem:Block2points}, the point of intersection is
%%   as well.

%% In the second case, it is still true that $B_P$ contains no fixed
%% points of $\tau^{2n}$ and at most one point from each $p$-cycle.
%% Therefore, if $|B_P| > 1$ then $B_P$ contains two points on distinct
%% lines and hence (by Proposition \ref{prop:trans-2} and Lemma  \ref{lem:Block2points}) every $\F_p$-point on a
%% smooth rational curve in characteristic $0$.  As in the first case it
%% follows that $B_P = \tK(\F_p)$, a contradiction, so $|B_P| = 1$.
%% All parts of a partition preserved by the action have the same size,
%% so $\mathcal B$ is discrete.
%% \end{proof}

The main result of this paper naturally breaks up into two
complementary parts.  Here we prove the first part of our main result
(Theorem~\ref{thm:alternating-kummer}).
%We will now prove the main result of this paper, which we restate
%below:
\begin{thm}\label{thm:alternating-kummer-deux} Let $E$ be an elliptic curve
over $\F_p$ with full level-$2$ structure,
where $p>2$, and let $K = \Kum(E \times E)$ be the Kummer 
surface of $E \times E$.  Let $\AI$ be the image of $\Azero$ in $\S_{K(\F_p)}$.  If $\AI$ is
transitive, then it contains the alternating group.
\end{thm}

%% \begin{thm}\label{thm:alternating-kummer-deux} Let $E$ be an elliptic curve
%% over $\F_p$ with full level-$2$ structure,
%% where $(p>2)$, and let $K = \Kum(E \times E)$ be the Kummer 
%% surface of $E \times E$.  Suppose that the action of the automorphism group 
%% of $K$ on the $\F_p$-rational points of $K$ is transitive.  Then the image
%% of $\Aut K$ in $\S_{K(\F_p)}$ contains the alternating group.
%% \end{thm}

The proof will be an application of a theorem of Praeger, which we
restate below for the convenience of the reader.  Note that Praeger's
result predates and does not depend on the classification of finite
simple groups.

\begin{thm}[{\cite[Theorem~A]{p}}]\label{thm:praeger}
Let $G$ be a primitive permutation group on a set of order $n$.
Suppose that $G$ has an element $\sigma$ of prime order $p$ with $q$ cycles of length $p$ and $f = n-qp$ fixed points.  Then one of the following holds:
\begin{enumerate}
\item $G \in \{\A_n, \S_n\}$;
 \item $G \in \{\A_c, \S_c\}$, where $n = c(c-1)/2$, $c = q+(p+1)/2$, and the action of $\A_c$ or $\S_c$ is that on the set of subsets of $\{1,\dots,c\}$ of order $2$;
 \item $n = 9, p = 3, q = 2$, and there are two possibilities for $G$;
 \item $f \le 5q/2 - 2$.
\end{enumerate}
\end{thm}

\begin{proof}(of Theorem \ref{thm:alternating-kummer-deux})
We apply Theorem \ref{thm:praeger} where $G$ is the image of $\Azero$
in the symmetric group action on $\tK(\F_p)$, $n=\#(\tK(\F_p)$ and
$\sigma=\tau^{2n}$.
Given these choices, $q = 16$ (cf. Lemma \ref{lem:good-power-of-tau}), so we are certainly not in case $3$.  In case $2$ we
would have $n = (p^2 + 64p + 1023)/8$.  However, $n$ is the number of
$\F_p$-points of $\tK$, which is at least $p^2+18p+1$, and this is
larger for $p> 7$.  In case $4$ we find $f \le 5q/2-2 = 38$, so $n \le
16p+38$, which again contradicts the estimate $n > p^2+18p+1$ for
$p > 5$.  Hence we are in case 1 as desired.  A simple computation (for
all possible surfaces) 
shows that the result holds for $p \le 7$.
\end{proof}

\begin{rmk} In the small cases $p = 3, 5$ where we have to compute directly,
  we use the method of Remark \ref{rem:how-to-find}.  That is, we apply a map
  to a point (possibly a point infinitely near to a node)
  by finding an elliptic curve smooth at that point and mapping an infinitely
  near point on that elliptic curve.
\end{rmk}

\subsubsection{Evenness of $\Azero$}\label{subsec:evenness-azero}

Recall that we refer to a permutation or a permutation group as {\em even}
if it is contained in the alternating group.
Here we complete the proof of our main result 
(Theorem~\ref{thm:alternating-kummer}) by proving its
complementary part:
\begin{thm}\label{thm:alternating-kummer-trois}
Let $E$ be an elliptic curve
over $\F_p$ with full level-$2$ structure,
where $p>2$, and let $K = \Kum(E \times E)$ be the Kummer 
surface of $E \times E$.
Let $\AI$ be the image of $\Azero$ in $\S_{K(\F_p)}$.
%If $\AI$ is
%transitive, then it
Then $\AI$ is even.
\end{thm}

\begin{lem}\label{lem:eight-multiple} %We first note that
The number of $\F_p$-points of $\tK$ is a multiple of
$8$.
\end{lem}
\begin{proof}
  We assumed that $E$ has full level-$2$ structure, so
  $a_p \equiv 2, 0 \bmod 4$ for $p \equiv 1, 3 \bmod 4$, where $a_p$ is the
  trace of Frobenius for $E$ over $\F_p$.  The number of points is
  $p^2+18p+1-a_p^2$, and the claim now follows easily.
\end{proof}
\begin{cor}\label{cor:even-4n-fixed}
  An automorphism of order $2$ is even in its action on $\tK(\F_p)$
  if and only if the number of $\F_p$-rational fixed points is a
  multiple of $4$.
\end{cor}

This only needs the number of $\tK(\F_p)$-points to be a multiple
of $4$, but we will need the more precise statement of
Lemma~\ref{lem:eight-multiple} in the proof of Theorem
\ref{thm:alternating-kummer-trois}.

\begin{rmk}\label{rem:parity-of-2-torsion}
  A $2$-torsion translation on an elliptic K3 surface
(or more generally a symplectic involution on any K3 surface) 
has $8$ fixed points, not necessarily defined over the ground field
\cite[Lemma 3]{nikulin}.
The singular points of singular irreducible fibres are fixed;
in addition, there may be fixed points on reducible fibres.  On the
other hand, an involution that is not a $2$-torsion translation may
fix various curves and isolated points.  For example, a negation map fixes
the $0$-section and the $2$-torsion as well as some components and points of
reducible fibres.
\end{rmk}

\begin{proof} (of Theorem \ref{thm:alternating-kummer-trois})
We prove that each of the automorphisms 
\begin{equation}
% f -> afaba and b -> ab
\{a, ab, c, d, ae, abga, abha, abia, afabcad, dj, cdk, c\ell\}
\end{equation}
acts by an even permutation.  In Section~\ref{alternategens} we noted as a 
consequence of the results of \cite[Section 4]{kk}
that $a, b, \dots, \ell$ generate $\Aut \tK$, and that immediately implies
that these elements do as well.

We start by proving that $a$ gives an even permutation (recall that $a$
is the map sending the image of $(x,y)$ on $E \times E$ to the image
of $(y,x)$).  The fixed points are precisely those on the curves
$(x,x)$ and $(x,-x)$; these do not intersect and each has $p+1$ 
points, so we get a total of $2p+2$.  This is a multiple of $4$, as is
the number of points on the Kummer surface (Lemma~\ref{lem:eight-multiple}).

%Translations

The permutation given by $ab$ is of order $3$, so it is even.

Let us
show that $4$ divides the number of fixed points over $\F_p$ for each
of the automorphisms $c, ae, abga$.  The argument for $d$ is the same as
for $c$, and for $abha, abia$ it is the same as for $abga$,
since the monodromy of the Legendre family realizes all permutations
of $\{g,h,i\}$ and fixes $a, b$.

First, $c$ is a translation for the constant fibration induced by the second
projection $E \times E \to E$, which has no
singular irreducible fibres.  It permutes the reduced components of
each reducible fibre without fixed points, so the $8$ fixed points
consist of $2$ on each nonreduced component.  The automorphisms of the
surface permute these, so either the fixed points are defined over
$\F_p$ on all the bad fibres or on none of them.  Thus there are $8$
or $0$, both of which are multiples of $4$.

Then $ae$ is a translation for the fibration $(x:y:z:w) \to (z:w)$,
whose reducible fibres are of type $I_2, I_2, I_8, I_8$.
The $4$ singular irreducible fibres are located at $s = \pm \sqrt \la$,
$\pm \sqrt{1/\la}$.  Thus they come in two pairs defined over the same
field $F(\sqrt \lambda)$ and either $0$ or $4$ fixed points are defined over
the base field.  The other $4$ fixed points are the points of
intersection of the components of the $I_2$ fibres.
There is an automorphism of  $\tK$ that preserves the set of fibres of
this fibration and exchanges the two $I_2$ fibres, namely $aca$, or
in coordinates $(x:y:z:w) \to (y:x:w:z)$.
Thus the last $4$ fixed points are all defined over the same field.
Again the number of fixed points is a multiple of $4$.

%bga -> abga
Now, $abga$ is a translation for a fibration with reducible fibres
$I_{12}, I^*_0$: in particular, the map can be defined by
$$(x^2 + yz + yw + zw: (\la-1)xy + (\la-1)y^2 - \la xz - yz - \la z^2 + xw + \la yw - (\la-1) zw + w^2).$$
There are no fixed points on the $I_{12}$, because no component goes to
itself or to an adjacent component.  On the $I_0^*$, there are
four fixed points: two intersection points of the nonreduced component
with a reduced component, and two on reduced components.  The first
type of point is fixed, so the second is as well, since they are
single points on curves defined over the base field.  Again there are
$4$ singular irreducible fibres in two pairs, this time over
$F(\sqrt{\lambda(\lambda-1)})$, where $F$ is the ground field.
Thus there are $4$ or $8$ fixed points over $F$.

%Negations:

% f -> afaba in this paragraph (it was 
We now show that $afabcad$ gives an even permutation.  Note that
$afabcad$ acts nontrivially on the base of a fibration with two
reducible $I_8$ fibres, namely the one defined by
$$(xy+yz-xw-zw:y^2-z^2).$$ 
Such a fibration has $8$ singular irreducible fibres for generic $\la$, whose
singularities are the only fixed points; their field of definition is
$F(\sqrt{-1},\sqrt{\la},\sqrt{\la-1})$.  So every cycle of $afabcad$ is of
length $4$ except for those involving these fixed points.
Not all of these points can be fixed by $afabcad$, because the action
on the base cannot have $8$ fixed points without being trivial.
For generic $\lambda$, they are all Galois conjugate, so if one of them
is fixed, so must the others be.  The points remain distinct for
$\la \notin \{0,1,\infty\}$, so by continuity none of them can be fixed
for any such $\la$.

Now, note that $(afabcad)^2 = cacda$ is $2$-torsion translation for the
given fibration and hence fixes all $8$ singular points.  It follows
that $afabcad$ has $4$ cycles of length $2$ on these points.  Since they are
in a single Galois orbit, all or none are defined over $\F_p$, and hence
we have either $0$ or $4$ cycles of length $2$.  
From Lemma~\ref{lem:eight-multiple} we know that $8|\#K(\F_p)$; all other
cycles are $4$-cycles, so there must be an even number of them.

The argument needs to be modified slightly when $\la+1$ or $\la^2-\la+1$ is $0$;
in these cases, two pairs of bad fibres coalesce to form $I_2$ fibres.
However, one can still calculate that $afabcad$ fixes none of these points
and has $0$ or $4$ cycles of length $2$.  To conclude, we note that
$a, b, c, d$ are all even, so $f$ is even if and only if $afabcad$ is.

The automorphism $dj$ fixes the class of a fibration with two
reducible fibres of type $I_4^*$, acting with order $2$ on the base
and exchanging the two reducible fibres.  The fibration is defined by
$$(2xy^2 - xz^2 + z^3 + 2xyw + 2yzw + xw^2  + zw^2:  2y^2z + xz^2 - z^3 + 2xyw + 2yzw + xw^2 + zw^2).$$

The two most obvious sections are the curves above the nodes at
$(-1:1:1:1)$ and $(1:1:-1:1)$; these are exchanged by $dj$, which acts
on the base by
$$(x:y) \to \left(\frac{\la+1}{\la-1}x-y:x-\frac{\la+1}{\la-1}y\right).$$
It acts as the identity on
the two fixed fibres, which are defined over the extension obtained by
adjoining $\sqrt \la$ to the ground field.  (These are located where one
of the generators of the Mordell-Weil group modulo torsion meets the zero
section.  There are no degeneracies to worry about, because the discriminant
of the polynomial whose roots are the fixed points is $16s/(s-1)^2$, which
is well-behaved at all $s$ for which the elliptic curve is smooth.)
If the two fixed fibres are not defined over the ground
field, clearly there are no fixed points.  
These fibres, unlike the general fibre, have full level-$2$ structure over
their field of definition, so the
number of fixed points over the field of definition is a multiple of $4$.
As before it follows that the permutation induced by $dj$ is even.
The arguments for $cdk$ and $c\ell$ are
virtually identical to this one.
\end{proof}
\section{Multiple Kummer surfaces}

One of the standard conjectures on Kummer  surfaces over number fields
is that every collection of finitely many local points can be
arbitrarily well approximated by a global rational point.  Analogously, one
might expect that every collection of finitely many
${\mathcal  O}_{K,\mathfrak p}$-points on a log K3 surface can be
arbitrarily well approximated by a global integral point.

As a step toward this result for the Markoff surface,
the authors of \cite{mpc} study the action of $\Gamma$ modulo several
primes at once, proving the following result among others:

\begin{thm} (\cite[Theorem 1.9]{mpc}) Let $n = p_1 \cdots p_k$ be a
product of distinct primes, and for $1 \le i \le k$ let $Q_{p_i}$ be the
permutation group induced by the action of the automorphism group of the
Markoff surface on points mod $p_i$ up to sign.  Suppose that all 
$Q_{p_i}$ are primitive permutation groups.  Then the induced action of
the automorphism group on $k$-tuples of points mod $p_i$ up to independent
sign changes is transitive.
\end{thm}

We prove a result analogous to theirs for the action of $\Azero$ on products
of Kummer surfaces.

\begin{prop}\label{thm:multiplekummers}
  Let $K_1, \dots, K_n$ be the Kummer surfaces of $E_i \times E_i$, where $E_i$ is
  an elliptic curve over a prime finite field $\F_{p_i}$
  % with $j$-invariant not in $\{0,1728,-3375\}$,
  and suppose that the
  action on each $K_i$ is transitive (hence alternating) and that, for
  all $i \ne j$, either $p_i \ne p_j$ or $K_i, K_j$ do not have the
  same number of $\F_{p_i}$-points.  Then the action of $\Aut K$ on the set of
  $n$-tuples $(P_i)$, where $P_i$ is an $\F_{p_i}$-point of $K_i$, is
  transitive, and indeed contains the product of the corresponding
  alternating groups.
\end{prop}

\begin{proof}Recall {\em Goursat's lemma} \cite[(4.19)]{suz}:
  let $G, G'$ be groups and let $H \subseteq G \times G'$
  be a subgroup for which both projections are
surjective.  Let $N, N'$ be the kernels of the projections of $H$ to
$G', G$ respectively.  Then $G/N \cong G'/N'$, with the image of $H$
being the graph of the isomorphism.

In particular, if $G, G'$ are nonisomorphic simple groups, then $N =
G, N' = G'$, and so $H = G \times G'$.  More generally, this follows
if the composition series for $G, G'$ have no simple groups in common,
and by a simple induction it follows that if the $G_i$ are pairwise
nonisomorphic simple groups and $H \subseteq \prod_{i=1}^n G_i$ is a
subgroup for which all projections are surjective then $H = \prod
G_i$.  So we are done, unless two of the $K_i$ (say $K_1, K_2$) have
the same number of points.

In this case, we have two possibilities in Goursat's lemma: either $N
= G, N' = G'$ as above, or $N, N'$ are both trivial.  Establishing
that the first of these holds requires exhibiting an element of $\Aut
K$ that acts trivially on $K_1$ but not $K_2$ (or vice versa).
Indeed, the element $\tau^2$ above has order $p_i$ on $K_i$, and our
assumption is that $p_1 \ne p_2$, so $\tau^{2p_1}$ fits the bill.
\end{proof}

\begin{rmk} Numerical computations suggest that Proposition
  \ref{thm:multiplekummers} holds even when some of the
$K_i$ are defined over the same field, as long as the elliptic curves
  are distinct.
\end{rmk}

\section{Appendix: Equations of fibrations used to define $k$ and
  $\ell$}\label{app:eqns}
%\adam{This should look like a section title: I don't know why} 
%\verb+\appendix+ \adam{doesn't work that way.}
In this appendix we present the equations for the elliptic fibrations used to
define the fibrations $k, \ell$, in analogy with those for $j$ given
in equation~(\ref{eqn:j}).

For $k$ we may define the
fibration by
\begin{align}\label{eqn:k}
\nonumber &(xy^3 + \frac{\la-2}{\la-1}xy^2z - \frac{1}{\la-1}y^3z - \frac{2\la}{\la-1}xyz^2 - 
\frac{2}{\la-1}y^2z^2 + 2x^2yw + \\
\nonumber &\quad \frac{2\la-1}{\la-1}xy^2w + 
\frac{\la}{\la-1}y^3w + \frac{2\la}{\la-1}x^2zw + \frac{6\la-4}{\la-1}xyzw + 
y^2zw + 2yz^2w - \\
\nonumber &\quad \frac{2}{\la-1}x^2w^2 - xyw^2 + 
\frac{2\la}{\la-1}y^2w^2 + \frac{\la-2}{\la-1}xzw^2 - \frac{2\la-1}{\la-1}yzw^2 + \\
&\quad \frac{1}{\la-1}xw^3 + \frac{\la}{\la-1}yw^3 - zw^3:\\
\nonumber &x^2yz + xy^2z + xyz^2 + y^2z^2 -  x^2yw - xy^2w + x^2zw - y^2zw + \\
\nonumber &\quad xz^2w + yz^2w - x^2w^2 - xyw^2 -  xzw^2 - yzw^2.)
\end{align}

For $\ell$, we use the fibration with equations 
\begin{align}\label{eqn:l}
\nonumber &(y^4 - \frac{2\la}{\la-1}xy^2z - \frac{2\la}{\la-1}y^3z -
\frac{2\la}{\la-1}xyz^2 - \frac{\la+1}{\la-1}y^2z^2 -
\frac{2\la}{\la-1}xy^2w - \frac{2\la}{\la-1}y^3w - \\
\nonumber &\quad \frac{4\la}{\la-1}xyzw - \frac{2\la+2}{\la-1}y^2zw -
\frac{2\la}{\la-1}xyw^2 - \frac{3\la-1}{\la-1}y^2w^2 +
\frac{2\la}{\la-1}xzw^2 + \\
&\quad \frac{2\la}{\la-1}yzw^2 - z^2w^2 +
\frac{2\la}{\la-1}xw^3 + \frac{2\la}{\la-1}yw^3 - 2zw^3: \\
\nonumber & xy^2z + y^3z + xy^2w + y^3w - xzw^2 - yzw^2 - xw^3 - yw^3).
\end{align}

\begin{thebibliography}{9}
\bibitem{acl} M. Artebani, C. Correa, A. Laface, {\em Cox rings of K3 surfaces of Picard number three}, J. Algebra 565, 598--626 (2021).

\bibitem{babai} L. Babai, {\em The probability of generating the symmetric group}, J. Combin. Theory Ser. A 52 (1), 148--153 (1989).

\bibitem{baragar} A. Baragar, {\em On the unicity conjecture for Markoff numbers}, Canad. Math. Bull. 39 (1), 3--9 (1996).

\bibitem{borcherds} R. Borcherds, {\em Automorphism groups of Lorentzian lattices}, J. Algebra 111, 133--153 (1987).

\bibitem{bgs} J. Bourgain, A. Gamburd, and P. Sarnak. {\em Markoff triples and strong approximation.}  C. R. Math. Acad. Sci. Paris 354 (2), 131--135. \url{arXiv:1607.01530}. 

\bibitem{cf} J. W. S. Cassels, E. V. Flynn, {\em Prolegomena to a middlebrow arithmetic of curves of genus 2.}  LMS Lecture Note Series 230.  Cambridge University Press, 1996.

\bibitem{frob} G. Frobenius, {\em \"Uber die Markoffschen Zahlen}, S. B. Preuss. Akad. Wiss, 458--487 (1913).

\bibitem{h-t} A. Harder, A. Thompson, {\em The geometry and moduli of K3 surfaces}.  In {\em Calabi-Yau varieties: arithmetic, geometry, and physics} (R. Laza, M. Sch\"utt, N. Yui, eds).  Fields Institute Monographs 34.  Springer, 2016.  \url{arXiv:1501.04049}.

\bibitem{isaacs} I. M. Isaacs, {\em Finite group theory.} Graduate Studies in Mathematics 92. AMS, 2008.

\bibitem{kk} J. Keum, S. Kondo, {\em The automorphism groups of Kummer surfaces associated with the product of two elliptic curves}, Trans. AMS 353 (4), 1469--1487 (2001).

\bibitem{kondo} S. Kondo, {\em The automorphism group of a generic Jacobian Kummer surface}, J. Algebraic Geometry 7, 589--609 (1998).

\bibitem{matsuki} K. Matsuki, {\em Introduction to the Mori program.}  Universitext.  Springer, 2002.

\bibitem{mpc} C. Meiri and D. Puder, with an appendix by D. Carmon.  {\em The Markoff group of transformations in prime and composite moduli}, Duke Math. J. 167 (14), 2679--2720 (2018).  \url{arXiv:1702.08358}.

\bibitem{nikulin} V. V. Nikulin, {\em On Kummer surfaces}, Izvestia Akad. Nauk SSSR 39, 278--293 (1975).

%\bibitem{nikulin2} V. V. Nikulin, {\em onElliptic fibrations on K3 surfaces},
%Proc. Edinburgh Math. Soc. 57 (1), 253--267 (2013).  \url{arXiv:1010.3904}.

\bibitem{od} E. O'Dorney, {\em Large orbits on Markoff-type K3 surfaces over finite fields}, \url{arXiv:2209.10436}.

\bibitem{p} C. Praeger,  {\em On elements of prime order in primitive permutation groups},  J. Algebra 60, 126--157 (1979).

\bibitem{shimada} I. Shimada, {\em An algorithm to compute automorphism groups of K3 surfaces and an application to singular K3 surfaces}, Int. Math. Res. Not. IMRN 2015 (22), 11961--12014 (2015).

\bibitem{shimada-web} I. Shimada, {\em Automorphism groups of complex elliptic K3 surfaces with Picard number 3}, \url{www.math.sci.hiroshima-u.ac.jp/shimada/preprints/AlgoAutK3/rho3examples/rho3examples.pdf}.

\bibitem{splag} J. H. Conway, N. J. A. Sloane, {\em Sphere packings, lattices, and groups} (third edition).  Grundlehren der mathematischen Wissenschaften 290.  Springer, 1999.

\bibitem{sterk} H. Sterk, {\em Finiteness results for algebraic K3 surfaces}, Math. Z. 189, 507--513, 1985.

\bibitem{suz} M. Suzuki, {\em Group theory I}.  Grundlehren der mathematischen Wissenschaften 247.  Springer, 1982.

\bibitem{tits} J. Tits, {\em Free subgroups in linear groups}, J. Algebra 20 (2), 250--270, 1972.

\bibitem{ujikawa} M. Ujikawa, {\em The automorphism group of the singular K3 surface of discriminant 7}, Comment. Math. Univ. St. Pauli 62 (1), 11--29 (2013).
\end{thebibliography}
\end{document}